\documentclass[11pt]{article}
\usepackage{amsmath,amssymb,mathrsfs,amsfonts}
\usepackage{hyperref}
\usepackage{sectsty}
\input xy
\xyoption{all}

\subsectionfont{\mdseries\itshape}

\subsubsectionfont{\mdseries\itshape}
\newtheorem{Thm}{Theorem}[section]
\newtheorem{Prop}[Thm]{Proposition}
\newtheorem{Lemma}[Thm]{Lemma}
\newtheorem{Cor}[Thm]{Corollary}

\newcommand{\pf}{\noindent{\em Proof.}\ }
\newcommand{\qed}{\hfill $\Box$}

\newcommand{\noin}{\noindent}

\newcommand{\ZZ}{\mathbb{Z}}
\newcommand{\QQ}{\mathbb{Q}}

\newcommand{\CC}{\mathbb{C}}
\newcommand{\NN}{\mathbb{N}}
\newcommand{\GG}{\mathbb{G}}
\renewcommand{\SS}{\mathbb{S}}
\newcommand{\TT}{\mathbb{T}}
\newcommand{\tensor}{\otimes}
\newcommand{\dc}{\cap}
\newcommand{\rc}{\subset}

\newcommand{\Spec}{\operatorname{Spec}}
\newcommand{\Spf}{\operatorname{Spf}}
\newcommand{\Fil}{\operatorname{Fil}}
\newcommand{\Gr}{\operatorname{Gr}}
\newcommand{\Art}{\mathfrak{Art}}

\newcommand{\OO}{{\cal O}}
\newcommand{\HH}{{\cal H}}
\newcommand{\UU}{{\cal U}}
\newcommand{\N}{{\cal N}}
\newcommand{\F}{{\cal F}}
\newcommand{\LL}{{\cal L}}
\newcommand{\MM}{{\cal M}}

\newcommand{\om}{\omega}
\newcommand{\oo}{\infty}

\newcommand{\lam}{\lambda}
\newcommand{\tha}{\theta}
\newcommand{\nb}{\nabla}
\newcommand{\pt}{\partial}
\newcommand{\eps}{\varepsilon}
\newcommand{\Hom}{\operatorname{Hom}}

\newcommand{\pd}[2]{#1^{[#2]}}
\newcommand{\ag}[1]{\langle #1 \rangle}
\newcommand{\dsb}[1]{[\hspace{-1.5pt}[ #1 ]\hspace{-1.5pt}]}

\newcommand{\unit}[1]{#1^{\times}}

\newcommand{\nnunit}[2]{\NN^{#1}\oplus\unit{#2}}
\newcommand{\logpt}[2]{{#1}_{(#2)}}
\newcommand{\tr}{\operatorname{tr}}

\newcommand{\Teich}{Teichm\"uller }
\renewcommand{\wp}{\frak{p}}

\begin{document}

\title{On Ordinary Crystals with Logarithmic Poles\footnote{
This work was partially supported by
the Golden-Jade fellowship of Kenda Foundation,
the NCTS (Main Office) and the NSC, Taiwan.}}
\author{\sc Jeng-Daw Yu\footnote{
Department of Mathematics,
National Taiwan University,
Taipei, Taiwan;
{\tt jdyu@math.ntu.edu.tw}.}}
\maketitle

\begin{abstract}
We derive some local properties of abstract crystals with logarithmic poles
over a smooth base in positive characteristic
and obtain the existence of the canonical coordinates
of certain ordinary crystals.
We then apply the results to deduce an integral property
of the coefficients of the so-called mirror maps.
\end{abstract}
%=========================================================
%=========================================================

Let $k$ be an algebraically closed field of characteristic $p > 0$.
The first aim of this paper is
to investigate the abstract formalism of crystals
over the ring $k\dsb{t_1, \cdots, t_m}$ of formal power series over $k$
with logarithmic poles (along some coordinate hyperplanes).
We then derive some of the basic properties of
the so-called ordinary crystals.
(Precise definitions are given in the text.)
Along the discussions,
the definition and the existence
of the ($p$-adic) canonical coordinates for an ordinary crystal
over $\Spec k\dsb{t_i}_{i=1}^m$ are obtained
(see Thm \ref{cancoord}).
The existence of these canonical coordinates,
which in some sense give a structural description
of the ordinary crystal,
is a generalization of a theorem of Deligne (\cite[\S 1.4]{Del_CC})
in the case without log structures.
Indeed the development of these results
is parallel to that of \cite{Del_CC}.
\smallskip

As an application,
we derive that the coefficients of the mirror map $\tilde{q}$
of a nice family of Calabi-Yau varieties over a number field
are $p$-adic integers for almost all prime numbers $p$.
(See \S \ref{Sect:App} for the precise statement.)
The integral property of the coefficients of $\tilde{q}$
has been observed in the 1990s.
In his letter to Morrison \cite{Del_inf} in 1993,
Deligne asked the question
whether the coefficients of $\tilde{q}$
can be written algebraically
in terms of the underlying family of Calabi-Yau varieties.
We do not know if this is true
or even how to formulate the statement properly.
On the other hand,
the integrality of $\tilde{q}$ was also observed independently
in \cite{HKTY} around the same time.
It has been proved in the first time in the paper \cite{LY}
for certain hypergeometric cases
using Dwork's criterion on the integrality of a $p$-adic formal power series.
In fact,
one of the main motivations of this work is the idea that
the mirror map $\tilde{q}$ (defined via the variation of Hodge structure)
should equal the canonical coordinates $q$
(in the sense of Serre-Tate local coordinates studied here),
and this relation should be revealed via the associated crystal.
It turns out that the two $q$'s are not exactly the same
but they differ by a multiplicative $p$-adic unit
(see Prop \ref{Prop:int_qc}).
However this suffices to deduce the integrality for $\tilde{q}$
(and also for the instanton numbers).
Since the mirror maps are defined at the place
where the fiber of the family is degenerate,
one introduces the log structures
to the family and to the crystal
and this leads to the results in this study.
\smallskip

The method of using Dwork's criterion
has been generalized in \cite{KR}
to the mirror maps obtained from certain differential equations.
This method has the advantage
that it does not put constrains on the prime $p$
and hence one obtains that the coefficients of the mirror maps
are indeed integers.
However to apply this method,
one has to know explicitly the coefficients
of the involved power series solutions
to the differential equation.
On the other hand,
the use of ordinary crystals
in the study of the mirror maps and the instanton numbers
has already appeared in the work \cite{JS}
in the non-singular case
and has been outlined in \cite{KSV} in the case with log poles.
The case of rank 4 and weight 3 has been treated in \cite{V}
as a follow-up of \cite{KSV}.
There to get the comparison between the two $q$'s mentioned above,
the author uses Voevodsky's category of mixed motives
and the theory of 1-motives
to extract information from the weight $\leq 1$ part of the rank four crystal,
which behaves like a variation of 1-motives.
However if the rank is larger,
this method does not apply.
\smallskip

The paper is organized as follows.
In the setup \S \ref{Sect:Kato},
we give the definitions of various crystals
over the scheme $\Spec k\dsb{t_i}_{i=1}^m$
equipped with the logarithmic structure
attached to the union $(\prod_{i=1}^r t_i)$
of the hyperplanes $\{(t_i)\}_{i=1}^r$.
In \S \ref{Sect:Local},
we focus on ordinary crystals
and derive the canonical coordinates of such a crystal (Thm \ref{cancoord}),
parallel to Deligne's theory in \cite{Del_CC}
for crystals with the trivial log structure.
\S \ref{Sect:App} consists of an application
concerning the integrality of the mirror map (Prop \ref{Prop:int_qc})
and of the instanton numbers (Thm \ref{Thm:kappa=kappa'})
of a certain family of Calabi-Yau varieties.
Finally in the appendix,
we briefly discuss the group structure\footnote{Thanks
are due to C.-L.~Chai for directing my attention to this question.}
on the underlying parameter space $S$
associated with a certain ordinary crystal
and the periods associated to the $W$-points $S(W)$.\medskip

\noin{\it Notations and conventions.}
Throughout this paper,
we let
$\NN =$ the additive monoid of non-negative integers
(with the unit 0).
We use the convention $0^0 = 1$, the multiplicative unit.

We fix an algebraically closed field $k$ of characteristic $p > 0$.
Let $W =$ the ring of Witt vectors with coefficients in $k$
and $K =$ the field of fractions of $W$.
Denote by $\sigma$ the absolute Frobenius
on $k, W$ and $K$.
Our universal base for the crystalline cohomology theory
of schemes over $k$
is the scheme $\Spec W$ with the trivial logarithmic structure $\unit{W}$.

\section{Crystals with logarithmic poles}\label{Sect:Kato}

In this section,
we give the definitions of various crystals
over $\Spec k\dsb{t_i}_{i=1}^m$ and over $\Spec k$ with logarithmic poles
and derive some relations among them.
For the corresponding results
of classical (i.e.~without log poles) crystals,
see \cite{Del_CC,Katz_Dwork}.
For definitions and more details regarding logarithmic structures,
see \cite{Kato_log}.

\subsection{The logarithmic schemes $S$ and $S_0$}

Throughout the paper\footnote{
The parameter $r$ indicates the number of degenerating directions
$\{t_i \in A\}_{i=1}^r$ on $S$,
which produce the non-trivial logarithmic structure.
When $r = 0$,
one recovers the definitions in the classical situation
where the involved log structure is the trivial one.},
we fix $r,m \in \NN$ with $r \leq m$.
Let $t = \{t_i\}_{i=1}^m$ denote the $m$-tuple of variables $t_i$.
Let $A = W\dsb{t}$ be the formal power series ring and
\begin{eqnarray}\label{def:L}
\LL = \nnunit{r}{A} &=& \bigcup_{n \in \NN^r} t^n \unit{A} \\
	(n,f) &\leftrightarrow& t^nf
\end{eqnarray}
where for $n  = (n_i)_{i=1}^r \in \NN^r$,
we set $t^n = \prod_{i=1}^r t_i^{n_i}$
as usual.
The monoid $\LL$ is generated over $\unit{A}$
by $\{ t_i \}_{i=1}^r$.
The natural inclusion $\LL \to A = \Gamma(\OO,\Spec A)$ of monoids
then defines a log structure on $\Spec A$.
We will use the same symbol $\LL$
to indicate the sheaf of monoids
underlying this log structure;
the sheaf $\LL$ is coherent in the sense of \cite[(2.1)]{Kato_log}.
Throughout the paper, let
\[ S = (\Spec A, \LL) \]
denote this logarithmic scheme.

To describe the module of logarithmic differentials,
which is generated by the log derivatives $d\log\LL$ of the monoid $\LL$,
it is convenient to
consider the twist $t' = \{t_i'\}_{i=1}^m$ of $t$ giving by
\begin{equation}\label{tprime}
t_i' = \left\{\begin{array}{ll}
	t_i & \text{if $1\leq i\leq r$} \\
	1+t_i & \text{if $r+1\leq i\leq m$}. \end{array}\right.
\end{equation}
The $A$-module $\om = \om^1_{S/W}$
of logarithmic differential 1-forms of $S$ over $W$
is generated freely by
\[ d\log t_i' = \frac{dt_i'}{t_i'} \quad (1\leq i \leq m). \]
Consider the derivations
\begin{equation}\label{def:tha}
\delta_i = \frac{d}{dt'_i} \quad\text{and}\quad
	\tha_i = \frac{d}{d\log t'_i} = t'_i\frac{d}{d t'_i}
		\quad (1\leq i \leq m).
\end{equation}
The $\tha = \{\tha_i\}_{i=1}^m$ is dual to $d\log t'$.
We have the following relations
\begin{equation}\label{rel:delta-theta}
(t'_i)^n \delta_i^n = \prod_{j=0}^{n-1} (\tha_i - j)
	\quad\text{for all $n\in \NN$}.
\end{equation}

Finally denote by
\[ S_0 = (\Spec A_0, \LL_0) \quad\text{and}\quad \om_0 = \om^1_{S_0/k} \]
the reductions mod $p$ of $S$ and $\om$, respectively.
The $\om_0$ is the module of log differentials
of the log scheme $S_0$ of characteristic $p$.

\subsection{Crystals over $S_0$}

Let $\HH$ be a free $A$-module.
With the appearance of the log structure,
a connection
\[ \nb: \HH \to \om \tensor_A \HH, \]
is called {\it quasi-nilpotent} (cf.~\cite[Remark (6.3)]{Kato_log}) if
for each $1\leq i\leq m$,
the endomorphism
\[ \prod_{j=0}^{n-1} (\nb(\tha_i) - j) \]
on $\HH$ tends to zero $p$-adically
as $n$ goes to $\oo$.
(Here $\tha_i$ is defined in \eqref{def:tha}.)\medskip

\noin{\it Definition.}
A {\em crystal} $(\HH,\nb)$ over $S_0$ is a pair
of a free $A$-module $\HH$
of finite rank
equipped with an integrable connection
\[ \nb: \HH \to \om \tensor_A \HH, \]
which is quasi-nilpotent.
\medskip

We now introduce the Frobenius and the Hodge structures on a crystal.
The absolute Frobenius $\sigma(a) = a^p$ on $A_0$
induces an endomorphism on the sub-monoid $\LL_0 \subset A_0$.
We call this $\sigma$ the {\it absolute Frobenius}
on $S_0$ or on $(A_0,\LL_0)$.
Notice that
a lifting $\phi: (A,\LL) \to (A,\LL)$ of $\sigma$ on $(A_0,\LL_0)$
must have the form
\[ \phi: t'_i \mapsto (t'_i)^p f_i
	\quad (1\leq i\leq m) \]
for some $f_i \in \unit{W\dsb{t}}$ with $f_i \equiv 1 \pmod{p}$.
(Here $t_i'$ is defined in \eqref{tprime}.)\medskip

\noin{\it Definition.}
An {\em F-crystal} $(\HH, \nb, \F)$ over $S_0$
is a triple consists of a crystal $(\HH,\nb)$ over $S_0$
with a Frobenius structure $\F$,
i.e.~for any lifting $\phi$ on $S$
of the absolute Frobenius,
there is a horizontal $A$-linear map
\[ \F(\phi): \phi^*\HH \to \HH, \]
which becomes an isomorphism
after $\tensor_W K$.
Here $\phi^*\HH := A \tensor_{\phi} \HH$.
For $x\in \HH$,
we write $\phi^*x = 1\tensor x \in \phi^*\HH$.\medskip

If $\phi_1$ and $\phi_2$ are two liftings of $\sigma$,
the compatibility of $\nb$ and $\F$ implies that
there is a unique isomorphism $\chi(\phi_1,\phi_2)$
such that
the following diagram commutes
\begin{equation}\label{diag:F-chi}
\xymatrix{
\phi_1^*\HH \ar[rr]^{\F(\phi_1)}\ar[d]_{\chi(\phi_1,\phi_2)} && \HH \\
\phi_2^*\HH \ar[urr]_{\F(\phi_2)} }.
\end{equation}
Indeed, if we write
\[\begin{array}{cc}
\phi_1(t'_i) = (t'_i)^p f_i \\
\phi_2(t'_i) = (t'_i)^p g_i \end{array}
	\quad(f_i,g_i \in 1+pW\dsb{t}, 1\leq i\leq m), \]
the isomorphism $\chi(\phi_1,\phi_2)$ is given explicitly by
\begin{equation}\label{1.1.3.4}
\chi(\phi_1,\phi_2)\phi_1^*x = \phi_2^*x
	+ \sum_{|n| > 0} \pd{\left(\frac{f}{g} - 1\right)}{n} \cdot
		\phi_2^*\left\{\left(\prod_{1\leq i \leq m}^{0\leq j \leq n_i -1}
			(\nb(\tha_i) - j)\right)x\right\}
				\quad\text{for $x \in \HH$}.
\end{equation}
Here $f = (f_i)_{i=1}^m, g = (g_i)_{i=1}^m$
and the index of the summation
runs through $n = (n_i) \in \NN^m$
with $|n| := n_1+\cdots +n_m > 0$.
Notice that we have
$f_ig_i^{-1} \equiv 1 \mod{p}$,
and hence the divided powers $\pd{(\bullet)}{n}: pA \to pA$
in the above formula are meaningful.
\medskip

\noin{\it Definition.}
A {\it Hodge $F$-crystal} $(\HH,\nb,\F,\Fil^{\bullet})$ over $S_0$
consists of an $F$-crystal $(\HH,\nb,\F)$
together with a decreasing filtration $\Fil^i \rc \HH$ of
free and co-free $A$-submodules
indexed by $i\in \NN$
such that
\begin{enumerate}
\item $\Fil^0 = \HH$ and $\Fil^n = 0$ for $n$ sufficiently large,
\item $\nb \Fil^{i+1} \rc \om\tensor_A \Fil^i$ for all $i$, and
\item $\F(\phi)\phi^*\Fil^i \rc p^i\HH$
for any lifting $\phi$ of the Frobenius and any $i$.
\end{enumerate}
The {\em (Hodge) weight} of such an $\HH$
is the sum $\rho_1 + \rho_2$
where
\begin{eqnarray*}
\rho_1 &=& \max\, \{n \in \NN\, | \Fil^n = \HH\}; \\
\rho_2 &=& \max\, \{n \in \NN\, | \Fil^n \neq 0\}.
\end{eqnarray*}

\subsection{Crystals over logarithmic points}

Consider the logarithmic scheme
\[ \Spec \logpt{W}{r} = (\Spec W, \nnunit{r}{W}) \]
where the log structure is defined by the morphism of monoids
\begin{eqnarray*}
\nnunit{r}{W} &\to& W \\
(n,x) &\mapsto& 0^n\cdot x.
\end{eqnarray*}
The $W$-module $\logpt{\om}{r}$ of differential 1-forms
is generated freely by the basis $\{\lam_i\}_{i=1}^r$.
Here if we let $\{\alpha_i\}_{i=1}^r$ be the standard basis of $\NN^r$,
then $\lam_i = d\log(\alpha_i,1)$,
where $(\alpha_i,1) \in \nnunit{r}{W}$.
\smallskip

Let
\[ \Spec \logpt{k}{r} = (\Spec k, \nnunit{r}{k}) \]
be the reduction of $\Spec\logpt{W}{r} \mod{p}$.
On the sheaf level, they fit into the commutative diagram
\[\xymatrix{
\nnunit{r}{W} \ar[r]\ar[d]_{({\rm id,\, mod}\, p)} & W \ar[d]^{{\rm mod}\, p} \\
\nnunit{r}{k} \ar[r] & k. }\]
From the diagram,
we see that there is a unique lifting,
still call it $\sigma$, to $\logpt{W}{r}$
of the absolute Frobenius on $\logpt{k}{r}$ given by
\begin{eqnarray*}
\sigma: \nnunit{r}{W} &\to& \nnunit{r}{W} \\
(n,u) &\mapsto& (pn, u^\sigma).
\end{eqnarray*}
Consequently we have
$\lam_i^{\sigma} = p\lam_i$ for all $i$
on the module $\logpt{\om}{r}$.
We remark that
the log scheme $\Spec \logpt{W}{r}$ is the
{\em canonical lifting} of $\Spec\logpt{k}{r}$
in the sense of \cite[Def (3.1)]{HK}.
\medskip

\noin{\it Definition.}
A {\em crystal} $(\HH, \N_i)$ over $\logpt{k}{r}$
is a pair consisting of a free $W$-module $\HH$ of finite rank,
endowed with
a collection of pairwise commutative,
quasi-nilpotent $W$-linear endomorphisms
$\{\N_i\}_{i=1}^r$.
An {\em F-crystal} $(\HH, \N_i, \F)$ over $\logpt{k}{r}$
is a crystal $(\HH, \N_i)$
together with a $\sigma$-linear map $\F: \HH \to \HH$,
which becomes a bijection after $\tensor_W K$,
and satisfies
\begin{equation}\label{rel:N-F}
\N_i\F = p\F\N_i
	\quad\text{for all $1\leq i\leq r$}.
\end{equation}
One defines the notion of a {\em Hodge $F$-crystals}
$(\HH, \N_i, \F, \Fil^{\bullet})$ over $\logpt{k}{r}$
and its weight
in a similar way as in the previous subsection
(e.g., $\N_i\Fil^{j+1} \rc \Fil^j$, etc.).\medskip

Now we explain the relation between the above definition
and the crystals obtained via connections.
We let ${\cal C}$ denote the category of crystals over $\logpt{k}{r}$.

Let ${\cal C}'$ be the category of pairs $(\HH,\nb)$
consisting of a free $W$-module $\HH$
with a quasi-nilpotent, integrable connection
\[ \nb: \HH \to \logpt{\om}{r}\tensor_W \HH. \]
We have a natural functor
\begin{eqnarray*}
{\cal C} &\to& {\cal C}' \\
(\HH, \N_i) &\mapsto&
	\left(\HH, \nb(h):=\sum_{i=1}^r \lam_i\tensor\N_i(h)\right).
\end{eqnarray*}
Here the commutativity of $\N_i$
translates to the integrability of $\nb$.
The following assertion is clear.

\begin{Prop}\label{Prop:C=C'}
The functor $C\to C'$ above establishes an equivalence of categories.
\end{Prop}

\subsection{The \Teich representative}

Fix an $l \in \NN$.
Any morphism $e_0: (A_0, \LL_0) \to \logpt{k}{l}$
determines uniquely a pair of mophisms
$(f_0,g_0): \NN^r \to \nnunit{l}{k}$ of monoids
obtained from the map between log structures
\begin{eqnarray*}
\LL_0 = \nnunit{r}{A_0} &\to& \nnunit{l}{k} \\
\big(n,a(t)\big) &\mapsto& \big(f_0(n), g_0(n)\cdot a(0)\big)
\end{eqnarray*}
with $f_0(n) = 0$ only if $n = 0$.
Conversely
any such a pair $(f_0,g_0)$ determines a map $e_0$.

\begin{Lemma}\label{Lemma:Teich-Base}
Let $\phi$ be a lifting to $S$
of the Frobenius $\sigma$ on $S_0$.
For any $e_0: (A_0, \LL_0) \to \logpt{k}{l}$,
there exists a unique lifting
$e = e(\phi): (A, \LL) \to \logpt{W}{l}$ of $e_0$
such that the following diagram commutes
\begin{equation}\label{diag:Teich}
\xymatrix{
(A, \LL) \ar[r]^e\ar[d]_{\phi} & \logpt{W}{l} \ar[d]^{\sigma} \\
(A, \LL) \ar[r]^e & \logpt{W}{l}. }
\end{equation}
\end{Lemma}

\pf
Write $e = (e^{\flat},e^{\sharp})$,
where $e^{\flat}: A \to W$ is a homomorphism
and $e^{\sharp}: \LL \to \nnunit{l}{W}$ is compatible with $e^{\flat}$.
The images $e^{\flat}(t_i)$ for $i > r$
are uniquely determined by $\phi$ as in the classical case
(see \cite[1.1.2]{Katz_Dwork}).

Let $(f_0,g_0): \NN^r \to \nnunit{l}{k}$ be the pair
corresponding to $e_0$ as above.
For $1\leq i\leq r$,
let $\alpha_i \in \NN^r \rc \LL$
corresponding to $t_i$ via \eqref{def:L}.
In this case,
we must have $e^{\flat}(t_i) = 0$
by a diagram chasing in \eqref{diag:Teich}.
On the level of log structures,
we have
\[ e^{\sharp}(\alpha_i,1) = (f_0(\alpha_i), x) \]
for some $x \in \unit{W}$
with
\begin{equation}\label{eqn:Teich-cong}
x \equiv g_0(\alpha_i) \mod{p}.
\end{equation}
Write $\phi(t_i) = t_i^p\gamma$ for some $\gamma \in 1+pA$.
Then the commutativity of \eqref{diag:Teich}
implies that
\[ x^{\sigma} = x^p\cdot \gamma(0), \]
which together with \eqref{eqn:Teich-cong}
determine $x$ uniquely.
\qed\medskip

\noin{\it Definition.}
Given a lifting $\phi$ of the Frobenius on $S$,
the (unique) lifting $e$, obtained in Lemma \ref{Lemma:Teich-Base},
of the morphism $e_0: (A_0, \LL_0) \to \logpt{k}{l}$
is called the {\em \Teich lifting} of $e_0$ relative to $\phi$.\medskip

The lifting $e$ induces a map
$\om^1_{S/W} \to \om^1_{\logpt{W}{l}/W}$
between differentials.
By the correspondence in Prop \ref{Prop:C=C'},
the {\em pull-back} $e^*\HH = \HH \tensor_A W$
of a (resp.~Hodge) $F$-crystal $\HH$ over $S_0$
by the \Teich lifting,
defines a (resp.~Hodge) $F$-crystal over $\logpt{k}{l}$.
The following lemma shows that
the pull-back of an $F$-crystal is independent of the lifting $\phi$.

\begin{Lemma}
Let $\HH$ be an $F$-crystal over $S_0$
and $e_0: (A_0,\LL_0) \to \logpt{k}{l}$ a morphism.
Suppose $\phi_1$ and $\phi_2$ are two liftings of the Frobenius
and let $e(\phi_1)$ and $e(\phi_2)$ be the corresponding
\Teich liftings.
Then the connection on $\HH$ provides
a canonical identification
$e(\phi_1)^*\HH = e(\phi_2)^*\HH$
as $F$-crystals over $\logpt{k}{l}$.
\end{Lemma}

\pf (Cf.~\cite[1.4]{Katz_Dwork})
As in the proof of Lemma \ref{Lemma:Teich-Base},
we write
\[\begin{array}{rcll}
e(\phi_i)^{\sharp}(\alpha_j) &=& (f_0(\alpha_j), \gamma_j(\phi_i)) &
	(1 \leq j \leq r) \\
e(\phi_i)^{\flat}(1+t_j) &=& \gamma_j(\phi_i) & (r < j \leq m)
\end{array}\]
to represent $e(\phi_i)$ for $i = 1,2$.
Let $\gamma(\phi_i) = \{\gamma_j(\phi_i)\}_{j=1}^m$.
Then the connection provides a canonical identification
(cf.~\eqref{1.1.3.4})
\begin{eqnarray}\label{eqn:e^*H}
e(\phi_1)^*\HH &\cong& e(\phi_2)^*\HH \\ \nonumber
e(\phi_1)^*\xi &\mapsto& \sum_{n \in \NN^m}
	\pd{\left(\frac{\gamma(\phi_1)}{\gamma(\phi_2)} -1\right)}{n}
	e(\phi_2)^*
		\left\{\left(\prod_{1\leq i \leq m}^{0\leq j \leq n_i -1}
			(\nb(\tha_i) - j)\right)\xi\right\}
				\quad\text{for $\xi \in \HH$.}
\end{eqnarray}
Let $\{\N_\mu(\phi_i)\}_{\mu = 1}^l$ be the associated
quasi-nilpotent operators on $e(\phi_i)^*\HH$.
Via Prop \ref{Prop:C=C'},
they equal the pull-backs by $e(\phi_i)$
of the residues of $\nb$ on $\HH$.
Thus we have
\[ \N_\mu(\phi_i) \big(e(\phi_i)^*\xi\big)
	= \sum_{\nu=1}^r \big[f_0(\alpha_\nu)\big]_\mu\cdot
		e(\phi_i)^*\big(\nb(\tha_\nu)\xi\big)
			\quad(1\leq \mu \leq l) \]
where $[\bullet]_\mu$ denotes the $\mu$-th components
of elements in $\NN^l$.
One computes trivially that
$\N_\mu(\phi_1)(e(\phi_i)^*\xi) \mapsto \N_\mu(\phi_2)(e(\phi_i)^*\xi)$
under the identification \eqref{eqn:e^*H} above.

On the other hand,
the induced Frobenius structures on the two pull-backs
are the same
by applying $e^*$ to the diagram \eqref{diag:F-chi}.
This completes the proof.
\qed\medskip

Thus by the above lemma,
we see that on $e^*\HH$,
the Frobenius structure and the $\N_j$
depend only on $e_0$,
but not the choice of $\phi$.
We denote this $F$-crystal by $e_0^*\HH$.
We remark that the Hodge filtration on $e_0^*\HH$
do depend on $\phi$ (cf.~Prop \ref{Prop:period}),
but the Hodge polygon does not.

\section{Ordinary crystals}
\label{Sect:Local}

\subsection{The setup}

\noin{\it Definition.}
An $F$-crystal $(\HH, \nb, \F)$ over $S_0$
is called a {\em unit-root $F$-crystal}
if for some (and hence for all) $\phi$ of the liftings of the absolute Frobenius,
$\F(\phi)$ is an isomorphism.
One defines the notion of
a {\em unit-root $F$-crystal over $\logpt{k}{r}$}
in a similar way.\medskip

Note that if $\HH$ is a unit-root $F$-crystal over $\logpt{k}{r}$,
then $\N_j$ are all trivial.
Indeed the equation \eqref{rel:N-F}
shows that each $\N_j$ is divisible by any power of $p$
and hence is identically zero.

On the other hand,
if $\HH$ is a unit-root $F$-crystal over $S_0$,
then there exists one basis $\{e_i\}_{i=1}^r$ of $\HH$ over $A$
such that $\nb e_i = 0$ and $\F(\phi)\phi^* e_i = e_i$
for all liftings $\phi$ of the Frobenius.
Thus $\HH$ is indeed an $F$-crystal over $A_0$
(i.e.~without poles).
Indeed one first constructs a basis $p$-adically inductively
verifying $\F(\phi)\phi^* e_i = e_i$ for a fixed lifting $\phi$.
Then $\{e_i\}$ satisfies the desired properties.
The proof is identical to that of \cite[Prop 1.2.2]{Del_CC}
after replacing $\Omega^1_{A/W}$ by $\om$
(granted the fact that
the log structure is irrelevant).\medskip

\noin{\it Definition.}\footnote{
In \cite{Del_CC},
such a crystal in this definition is called
an {\em ordinary Hodge $F$-crystal}.
Here we shorten the terminology
for abbreviation.}
Let $e_0: S_0 \to \Spec \logpt{k}{r}$ be the augmentation.
A Hodge $F$-crystal $\HH$ over $S_0$ is called an
{\em ordinary crystal}
if the Newton and the Hodge polygons of the $F$-crystal
$e_0^*\HH$ over $\logpt{k}{r}$ coincide.
(Cf.~\cite[Prop 1.3.2]{Del_CC}.)

\begin{Prop}\label{Prop:UH}
Let $\HH$ be an ordinary crystal over $S_0$.
There exists a unique filtration of $\HH$ by sub-$F$-crystals
\begin{equation}\label{def:U}
0 \rc \UU_0 \rc \UU_1 \rc \cdots \rc \UU_i \rc \UU_{i+1} \rc \cdots
\end{equation}
which satisfies the following two properties:
\begin{enumerate}
\item
$\UU_i/\UU_{i+1}$ is of the form ${\cal V}_i(-i)$,
where ${\cal V}_i$ is a unit-root $F$-crystal and $(-i)$
denotes the Tate twist.
\item
We have the decompositions of the $A$-module
\[ \HH = \UU_i \oplus \Fil^{i+1} \quad\text{for all $i$} \]
and 
\begin{equation}\label{def:Hi}
\HH = \bigoplus_{i\in \NN} \HH^{(i)}
	\quad\text{where $\HH^{(i)} := \UU_i \dc \Fil^i$}.
\end{equation}
\end{enumerate}
In particular,
if $e_0: S_0 \to \Spec \logpt{k}{r}$ is the augmentation,
then $\N_j^{\rho+1} = 0$ on $e_0^*\HH$ for all $j$,
where $\rho =$ the weight of $\HH$.
\end{Prop}

\pf
The assertions (i) and (ii) are proved identically as in the classical case,
see \cite[\S 1.3]{Del_CC}.
Indeed one constructs the filtration $\UU_\bullet$ inductively
by using the Frobenius structure,
which itself does not see the log structure,
and reducing to the unit-root case
where the log structure is irrelevant by the discussion above.

Since $\N_j$ acts trivially on the successive quotients
$e_0^*(\UU_i/\UU_{i+1})$,
we have $\N_j^{\rho+1} = 0$.
\qed

\subsection{Ordinary crystals of weight one}

Let $\HH$ be an ordinary crystal of weight one
over $S_0$.
Write $\HH = \UU \oplus \Fil^1$ as in Prop \ref{Prop:UH},
where $\UU$ is the unit-root part of $\HH$.
The following theorem describes the structure of $\HH$;
the proof will be given in the next subsection.

\begin{Thm}[Cf.~{\cite[Th 1.4.2]{Del_CC}}]\label{cancoord}
With assumptions as above, we have the following.
\begin{enumerate}
\item
There exist bases
$a = \{a_i\}_{i=1}^g$ and $b = \{b_j\}_{j=1}^h$
of the $A$-modules $\UU$ and $\Fil^1$, respectively,
verifying
\begin{equation}\label{spade1}
\begin{array}{rcll}
\nb a_i &=& 0 & (1\leq i\leq g) \\
\nb b_i &=& \sum_{j=1}^g \eta_{ij}\tensor a_j
	 & (1\leq i\leq h), \\
\end{array}
\end{equation}
for some closed
$\eta_{ij} \in \om^1_{S/W}$,
and for every lifting $\phi$ of the Frobenius,
\begin{equation}\label{spade2}
\begin{array}{rcll}
\F(\phi)\phi^*a_i &=& a_i & (1\leq i\leq g) \\
\F(\phi)\phi^*b_i &=& pb_i + p\sum_{j=1}^g u_{ij}(\phi)a_j
	& (1\leq i\leq h),
\end{array}
\end{equation}
for some $ u_{ij}(\phi) \in A$.
Moreover
$\eta_{ij}$ and $u_{ij}(\phi)$ satisfy
\begin{equation}\label{heart}
\phi^* \eta_{ij} = p\eta_{ij} + p\cdot du_{ij}(\phi)
	\quad\text{for each $\phi$}.
\end{equation}
\item
With the choice of $a, b$ above,
there exists a unique collection of elements
\[ \tau_{ij} = \tau^{\log}_{ij} + \tau'_{ij}
	\quad((1,1)\leq (i,j)\leq (g,h)) \]
with
\[ (\tau^{\log}_{ij} , \tau'_{ij}) \in
	\left(\bigoplus_{l=1}^r \ZZ_p \cdot \log t_l\right) \oplus K\dsb{t} \]
such that for all $i, j$ and $\phi$,
\begin{equation}\label{def:tau}
\eta_{ij} = d\tau_{ij} \quad\text{with $\tau'_{ij}(0) \in pW$}
\end{equation}
\begin{equation}\label{club}
\phi^* \tau_{ij} - p\tau_{ij} - pu_{ij}(\phi) = 0.
\end{equation}
\item
Suppose $p\neq 2$.
The power series
\begin{equation}\label{def:q'}
q_{ij}' = \exp(\tau'_{ij}) = \sum\frac{(\tau'_{ij})^n}{n!}
	\quad((1,1)\leq (i,j)\leq (g,h))
\end{equation}
are well-defined elements in $A$,
and verify $q_{ij}'(0) \equiv 1 \mod{pW}$.
\end{enumerate}
\end{Thm}

Keep the $\HH = \UU \oplus \Fil^1$ over $S_0$
as above.
The connection provides a linear map
via the composition
\[ \Fil^1 \overset{\nb}{\longrightarrow} \om \tensor \Fil^0
	\to \om \tensor (\Fil^0/\Fil^1)
	\overset{\sim}{\longrightarrow} \om \tensor\UU. \]
Let $\om^{\vee}_{S/W}$ denote the dual of the $A$-module $\om$.
We then obtain an $A$-linear map
\begin{equation}\label{dag}
\Gr\nb: \om^{\vee}_{S/W} \to \Hom_A (\Fil^1, \UU),
\end{equation}
which in some sense measures the deviation
of the crystal $\HH$ over $S_0$.
In a special case that will occur in the next section,
we have the following.

\begin{Cor}[Cf.~{\cite[Cor 1.4.7]{Del_CC}}]\label{Cor:wt1}
With the notations in Theorem \ref{cancoord},
suppose $p$ is odd, $g = 1$ and the map \eqref{dag} is an isomorphism.
Then we can further modify the basis $\{a,b\}$ of $\HH$ above
such that, in addition to the statements therein,
the following hold :
\begin{enumerate}
\item
With $q'_{1j}$ defined in \eqref{def:q'},
let
\begin{equation}\label{eqn:def:q}
q_j = \left\{\begin{array}{cl}
	t_jq'_{1j} & \text{if $1\leq j\leq r$} \\
	q'_{1j} & \text{if $r < j\leq m$}. \end{array}\right.
\end{equation}
The $W$-homomorphism
\begin{eqnarray*}
(W\dsb{x}, \LL_x) &\to& (A,\LL) \\
x_j &\mapsto& \left\{\begin{array}{cl}
	q_j & \text{if $1\leq j\leq r$} \\
	q_j - 1 & \text{if $r < j\leq m$} \end{array}\right.
\end{eqnarray*}
is an isomorphism of affine log schemes,
where $x = \{x_j\}_{j=1}^m$ and
\[ \LL_x = \bigcup_{n \in \NN^r} x^n\cdot \unit{W\dsb{x}}. \]
\item
Let $\phi$ be the lifting of the Frobenius
defined by $\phi(q_i) = q_i^p$.
We have
\[ \F(\phi)\phi^* b_i = pb_i
	\quad\text{for all $1\leq i\leq m$.} \]
\end{enumerate}
Furthermore,
the choice of $\{q_i\}_{i=1}^r$
depends only on the arrangement of $\{t_i\}_{i=1}^r$
modulo $\unit{(W\dsb{t_i}_{i=r+1}^m)}$ multiplicatively.
\end{Cor}

\pf
Since \eqref{dag} is a bijection,
we can change the bases $a$ and $\{b_j\}$
by multiplying by elements in $\unit{\ZZ_p}$ and ${\rm GL}_m(\ZZ_p)$,
respectively
such that $\tau_{1j}^{\log} = \log t_j$
for $1\leq j\leq r$.
Thus $\exp(\tau_{1j}) = t_j\exp(\tau_{1j}')$
are well-defined power series in $A$.
(Note that the modification does not destroy the relations
\eqref{spade1} and \eqref{spade2}.)
The statement (i) is now clear.

The assertion (ii) follows from \eqref{heart}
as in this case, $\phi^*\eta_{1j} = p\eta_{1j}$.
\qed

\subsection{The proof of Theorem \ref{cancoord}}
To prove the assertions in (i),
one reduces to the unit-root case,
where the log structure is irrelevant,
and thus can proceed identically
as in the classical case
(see \cite[Th 1.4.2(i)]{Del_CC}
and replace $\Omega^1_{A/W}$ by $\om$ there).
Before the proof of (ii),
we need the following lemma,
which is the log analogue of the Poincar\'e lemma.

\begin{Lemma}\label{Poincare}
Let $x = \{x_i\}_{i=1}^g$ and $y = \{y_j\}_{j=1}^h$
be distinct variables
and $\eta = \sum \alpha_i d\log x_i + \sum \beta_j dy_j$
be a closed differential form with
$\alpha_i, \beta_j \in W\dsb{x,y}$.
We have $\eta = d \tau$
for a unique
$\tau \in W \cdot \log x \oplus K\dsb{x,y}$
up to an additive constant.
Moreover
$\tau \equiv \sum \alpha_i\log x_i \mod{K\dsb{x,y}}$.
\end{Lemma}

\pf
By the assumption,
\begin{equation}\label{eta_closed}
\begin{split}
0 = d \eta = \sum_{i,j} \left( x_i \frac{\pt \alpha_j}{\pt x_i}
			- x_j \frac{\pt \alpha_i}{\pt x_j} \right)
				\frac{dx_i}{x_i} \wedge \frac{dx_j}{x_j}
		+ \sum_{j,l} \left( x_j \frac{\pt \beta_l}{\pt x_j}
			- \frac{\pt \alpha_j}{\pt y_l} \right)
				\frac{dx_i}{x_i} \wedge dy_l \\
		+ \sum_{l,n} \left( \frac{\pt \beta_n}{\pt y_l}
			- \frac{\pt \beta_l}{\pt y_n} \right)
				dy_l \wedge dy_n.
\end{split}\end{equation}
Pick any
\[ \tau_1 = \int \frac{\alpha_1}{x_1} dx_1
	\in W\cdot \log x_1 \oplus K\dsb{x,y}. \]
Notice that this is possible
for the equation \eqref{eta_closed} above implies that
$\alpha_1 - \alpha_1(0)$ is divisible by $x_1$
as an element in $W\dsb{x,y}$.
Let $\eta_1 = \eta - d \tau_1$.
Then $\eta_1$ remains closed and
\[ \eta_1 = \sum_{i=2}^g \left( \alpha_i - x_i \frac{\pt \tau_1}{\pt x_i}\right)
		\frac{dx_i}{x_i}
	+ \sum_{j=1}^h \left( \beta_j - \frac{\pt \tau_1}{\pt y_j}\right)
		dy_j \]
does not have the term involving $\frac{dx_1}{x_1}$.
Moreover, we have
\[ \frac{\pt}{\pt x_1}
		\left( \alpha_i - x_i \frac{\pt \tau_1}{\pt x_i}\right)
	= \frac{1}{x_1}\left( x_1 \frac{\pt \alpha_i}{\pt x_1}
		- x_i \frac{\pt \alpha_1}{\pt x_i} \right) = 0 \]
by \eqref{eta_closed} and similarly
\[ \frac{\pt}{\pt x_1}
		\left( \beta_j - \frac{\pt \tau_1}{\pt y_j}\right)
	= \frac{1}{x_1}\left( x_1 \frac{\pt \beta_j}{\pt x_1}
		- \frac{\pt \alpha_1}{\pt y_j} \right) = 0. \]
Thus $\eta_1$ is independent of $x_1$.
Now by induction and the classical Poincar\'e lemma
(i.e.~the case $g=0$),
the assertion follows.
\qed\medskip

Now back to the proof of Thm \ref{cancoord}.
\medskip

{\it Proof of \eqref{def:tau} ---}
An inspection of \eqref{heart} reveals that
the residues of $\eta_{ij}$ are in $\ZZ_p$.
Together with Lemma \ref{Poincare},
we then see that
$\eta_{ij} = d\tau_{ij}$ for some
$\tau_{ij} \in \ZZ_p \cdot \log t \oplus K\dsb{t}$.
We further write
\begin{eqnarray*}
\tau_{ij} &=& \tau^{\log}_{ij} + \tau'_{ij} \\
	&=& \left(\sum_{l=1}^r \tau^{(l)}_{ij} \cdot \log t_l\right) + \tau'_{ij}
		\quad (\tau^{(l)}_{ij} \in \ZZ_p)
\end{eqnarray*}
corresponding to the direct sum decomposition
as in the statement in (ii).

Recall the $t' = \{t'_l\}$ defined in \eqref{tprime}.
For any lifting $\phi$,
we have
\[ \phi(t'_l) = (t'_l)^p f_l(\phi) \quad (1\leq l\leq m) \]
for some
\[ f(\phi) = (f_l(\phi)) \in \big(1 + pW\dsb{t}\big)^m. \]
In particular
\[ \log f_l(\phi) := -\sum_{n=1}^\infty \frac{(1-f_l(\phi))^n}{n} \]
are well-defined elements in $pW\dsb{t}$.
Define $v_{ij}(\phi) \in W\dsb{t}$
by letting
\begin{equation}\label{def:v}
pv_{ij}(\phi) = pu_{ij}(\phi) - \sum_{l=1}^r \tau^{(l)}_{ij} \log f_l(\phi).
\end{equation}
By substituting \eqref{def:v} into the integration of \eqref{heart},
we see that
\[ \phi^*\tau'_{ij} - p\tau'_{ij} - pv_{ij}(\phi) \in K \]
is a constant and hence
\[ \phi^*\tau'_{ij} - p\tau'_{ij} - pv_{ij}(\phi)
	= \left(\phi^*\tau'_{ij}\right)(0) - p\tau'_{ij}(0) - pv_{ij}(\phi)(0). \]
Now let $\psi$ be the lifting given by sending $t$ to $t^p$
(and thus $v_{ij}(\psi) = u_{ij}(\psi)$).
We normalize $\tau_{ij}$ by requiring
\begin{equation}\label{bigstar}
\left(\psi^*\tau'_{ij}\right)(0) - p\tau'_{ij}(0) - pv_{ij}(\psi)(0) = 0.
\end{equation}
Since in this case,
$\left(\psi^*\tau'_{ij}\right)(0) = \tau'_{ij}(0)^{\sigma}$,
we have explicitly
\[ \tau'_{ij}(0) = \sum_{n=1}^\infty V^n v_{ij}(\psi)(0)
	\quad\text{where $Vx := px^{\sigma^{-1}}$.} \]
This shows that $\tau'_{ij}(0) \in pW$.
\medskip

{\it Proof of \eqref{club} ---}
It suffices to verify that for any lifting $\phi$ of the Frobenius,
the equality \eqref{bigstar} remains valid
after replacing $\psi$ by $\phi$.
We distinguish $\phi$ into two cases:

{\em Case 1. The augmentation map
$ e: \logpt{W}{r} \to S $
equals the \Teich representative relative to $\phi$
of the augmentation $e_0 = e \tensor_W k$.}
In this case,
\[ f_l(\phi)(0) = 1
	\quad\text{for all $1\leq l\leq m$.} \]
Thus $v_{ij}(\phi)(0) = u_{ij}(\phi)(0)$
and the two maps $\F(\phi)\phi^*$ and $\F(\psi)\psi^*$
induce the same Frobenius on $e^*\HH$.
Applying $e^*$ to the the second equation in \eqref{spade2},
we obtain
\begin{equation*}
p(e^*b_i) + p \sum_{j=1}^g u_{ij}(\phi)(0) (e^*a_j) = \F(e^*b_i)
	= p(e^*b_i) + p \sum_{j=1}^g u_{ij}(\psi)(0) (e^*a_j).
\end{equation*}
Hence $u_{ij}(\phi)(0) = u_{ij}(\psi)(0)$.

{\em Case 2. For general $\phi$.}
We have
\[ \phi(t_l') = (t'_l)^p f_l'(\phi) \cdot (1 + pw_l^{\sigma})
	\quad(1\leq l\leq m) \]
for some $w = (w_l) \in W^m$ and $f_l'(\phi)(0) = 1$.
Let $\phi_0$ be the Frobenius with $\phi_0(t'_l) = (t'_l)^p f'_l(\phi)$.
Then $\phi_0$ is compatible with the augmentation $e$.

Applying $\chi(\phi,\phi_0)$ in \eqref{1.1.3.4}
to the second equation for $\phi$ and $\phi_0$ in \eqref{spade2},
one obtains
\begin{equation}\label{star}
p\sum_{j=1}^g u_{ij}(\phi)a_j = p\sum_{j=1}^g u_{ij}(\phi_0)a_j \\
	+ \sum_{|n| > 0}
	\pd{(pw^{\sigma})}{n}\cdot \F(\phi_0)\phi_0^*
		\left\{\left(\prod_{1\leq l \leq m}^{0\leq j \leq n_l -1}
			(\nb(\tha_l) - j)\right)b_i\right\}.
\end{equation}
Using the second equation of \eqref{spade1}
and plugging \eqref{def:tau} in,
one gets
\[ \left(\prod_{1\leq l \leq m}^{0\leq j \leq n_l -1}
	(\nb(\tha_l) - j)\right)b_i
		= \sum_{j=1}^g (t')^n (\delta^n\tau_{ij})a_j
			\quad (1\leq i\leq h), \]
where $\delta_l = \frac{d}{dt'_l}$
(see Equation \eqref{rel:delta-theta}).
Thus, from \eqref{star} and the substitution \eqref{def:v},
we deduce
\[ pv_{ij}(\phi) =
	pv_{ij}(\phi_0) - \sum_{l=1}^r \tau^{(l)}_{ij} \log(1+ pw^{\sigma})
	+ \sum_{|n| > 0} \pd{(pw^{\sigma})}{n}
		\phi_0^*\left((t')^n\cdot \delta^n\tau_{ij}\right). \]
Therefore,
to prove \eqref{bigstar} for $\phi$,
it reduces to show that
\begin{equation*}
\begin{split}
\left(\phi^*\tau'_{ij}\right)(0) - p\tau'_{ij}(0) - pv_{ij}(\phi_0)(0) \\
	+ \sum_{l=1}^r \tau^{(l)}_{ij} \log(1 + pw_l^{\sigma})
	&- \sum_{|n| > 0} \pd{(pw^\sigma)}{n}
		\left((t')^n\cdot \delta^n\tau_{ij}\right)(0)^{\sigma} = 0.
\end{split}
\end{equation*}

For the last term, we have
\begin{eqnarray*}
\sum_{|n| > 0} \pd{(pw^\sigma)}{n}
	\left((t')^n\cdot \delta^n\tau_{ij}\right)(0)^\sigma
	&=& \sum_{l=1}^r (\tau^{(l)}_{ij})^\sigma \sum_{n_l > 0}
		(-1)^{n_l+1}\pd{(pw_l^\sigma)}{n_l} \\
	&=& \sum_{l=1}^r \tau^{(l)}_{ij} \log(1 + pw_l^{\sigma}),
\end{eqnarray*}
since $\tau^{(l)}_{ij} \in \ZZ_p$.
Thus it suffices to show that
\[ \left(\phi^*\tau'_{ij}\right)(0) - p\tau'_{ij}(0) - pv_{ij}(\phi_0)(0) = 0. \]

Notice that $\phi = \phi_0 \circ \alpha$,
where $\alpha$ is the $W$-endomorphism of $A$
defined by
\[ \alpha(t'_l) = t'_l(1 + pw_l')
	\quad\text{where $w_l' := \sigma^{-1}(w_l)$}. \]
Hence we have
\[ \left(\phi^*\tau'_{ij}\right)(0)
	= \left(\phi_0^*(\alpha^*\tau'_{ij})\right)(0)
	= \left(\alpha^*\tau'_{ij}\right)(0)^{\sigma}
	= \tau'_{ij}(0)^{\sigma}. \]
Granting the validity of \eqref{bigstar} for $\phi_0$,
we have
\[ \tau'_{ij}(0)^{\sigma} = p\tau'_{ij}(0) + pv_{ij}(\phi_0)(0). \]
Thus the equality \eqref{bigstar} holds for $\phi$.
\medskip

{\it The rest ---}
The proof of the properties of $q'_{ij}$ in (iii),
which is a consequence of \eqref{club} and $\tau'_{ij}(0) \in pW$
by a trick of Dwork,
is identical to that of \cite[Th 1.4.2(ii)]{Del_CC}.
The proof of the theorem is now complete.
\qed

\section{Ordinary crystals of Calabi-Yau type}
\label{Sect:App}

\subsection{The Calabi-Yau condition}

\noin{\it Definition.} A Hodge $F$-crystal $(\HH,\Fil^{\bullet})$
of weight $\rho$ over $S_0$
is called of {\em Calabi-Yau type}
if $\Fil^\rho$ is of rank 1 over $A$
and $\HH$ is equipped with a horizontal, perfect,
$(-1)^\rho$-symmetric pairing
of $F$-crystals
\[ \ag{\, ,}: \HH\times\HH \to A(-\rho) \]
such that
$\ag{\Fil^i,\Fil^{\rho+1-i}} = 0$
for all $i$.\medskip

Suppose $\HH$ is ordinary and of Calabi-Yau type.
Then it is easy to deduce from the definition that
$\ag{\UU_i, \UU_{\rho-1-i}} = 0$ for all $i$,
where $\UU_{\bullet}$ is the filtration in \eqref{def:U}.
For the following theorem,
cf.~\cite[Th 2.1.7]{Del_CC}, \cite[\S 14]{Del_inf}
and \cite[Thm 2.2]{JS}.

\begin{Thm}\label{Thm:weight=3}
Let $\HH$ be an ordinary crystal of Calabi-Yau type of weight three
over $S_0$.
Suppose $p \neq 2$
and the map \eqref{dag} is an isomorphism.
Then there exist a basis $\{u_0, u_1^{(i)}, u_2^{(i)}, u_3\}_{i=1}^m$
of $\HH$
%adapted to the decomposition $\HH = \bigoplus \HH_i$
with $\ag{u_0,u_3} = 1 = \ag{u_1^{(i)},u_2^{(i)}}$,
and elements $q_i \in A$
verifying the following properties:
\begin{enumerate}
\item
$A = W\dsb{q_i, q_j-1}_{1\leq i\leq r}^{r < j\leq m}$
with
\[ \left.\begin{array}{cl}
	(t_i^{-1}q_i)(0) & (1\leq i\leq r) \\
	q_j(0) & (r < j\leq m) \end{array}\right\}
		\equiv 1 \mod{pW}. \]
\item The connection is given by
\[\begin{array}{llll}
\nb u_0 &=& 0 \\
\nb u_1^{(i)} &=& d\log q_i \tensor u_0 \\
\nb u_2^{(i)} &=& \sum_j \beta_{ij} \tensor u_1^{(j)}
	&\text{for some $\beta_{ij} \in \om^1_{S/W}$} \\
\nb u_3 &=& -\sum_i d\log q_i \tensor u_2^{(i)} . \end{array} \]
\item
Write
\begin{equation}\label{eqn:def:kappa}
\beta_{ij} = \sum_l \kappa_{ijl} \cdot d\log q_l
	\quad(\kappa_{ijl} \in A).
\end{equation}
If $\phi: A \to A$ is the lifting of the Frobenius given by $\phi(q_i) = q_i^p$,
then
\[\begin{array}{lll}
\F(\phi)\phi^*u_0 &=& u_0 \\
\F(\phi)\phi^*u_1^{(i)} &=& pu_1^{(i)} \\
\F(\phi)\phi^*u_2^{(i)} &=& p^2(a_i u_0 + \sum_j b_{ij} u_1^{(j)} + u_2^{(i)}) \\
\F(\phi)\phi^*u_3 &=& p^3(cu_0 + \sum_i a_i u_1^{(i)} + u_3) \end{array} \]
for some $a_i, b_{ij}, c \in A$.
Let
$\pt_i = q_i\frac{\pt}{\pt q_i}$.
We have
\begin{equation}\label{abck}
\begin{array}{lll}
\pt_i c &=& -2a_i \\
\pt_j a_i &=& b_{ij} \\
\pt_l b_{ij} &=& \phi(\kappa_{ijl}) - \kappa_{ijl}.
\end{array}
\end{equation}
\end{enumerate}
\end{Thm}

\pf
Let $\UU_i, \HH_i$ be defined as in Prop \ref{Prop:UH}.
Then $\UU_2$, with the induced filtration from $\Fil^{\bullet}$,
is ordinary of weight 1.
Thus by Cor \ref{Cor:wt1},
there exist $u_0, u_1^{(i)}$ and $q_i$
satisfying the desired conditions
in (i), (ii) and (iii).

Since $\UU_0^{\perp} = \UU_2$ and $\UU_1^{\perp} = \UU_1$,
we can find a basis $\{u_i\}$ of $\HH_2$
such that
$\ag{u_1^{(i)}, u_i} = 1$.
Let
\[ u_i' = u_i - \sum_{j \neq i} \ag{u_1^{(j)}, u_i} u_j. \]
Then $\ag{u_1^{(j)},u_i'} = 0$ if $i\neq j$
and $\{ u_i'\}$ forms a basis of $\HH_2$.
Since the product $\ag{\,,}$ induces a perfect pairing
between $\HH_1$ and $\HH_2$
and $\ag{\HH_1,\HH_1} = 0 = \ag{\HH_2,\HH_2}$,
the elements $\ag{u_1^{(i)},u_i'}$ must be invertible.
Thus the set
\[ \left\{ u_2^{(i)} = \ag{u_1^{(i)},u_i'}^{-1}\cdot u_i' \right\} \]
forms a basis of $\HH_2$,
satisfying $\ag{u_1^{(i)},u_2^{(j)}} = \delta_{ij}$.

Take $\ag{u_3} = \HH_3$ with $\ag{u_0, u_3} = 1$.
By the transversality of the Hodge filtration
and that $\nb(\UU_{i+1}/\UU_i) =0$,
we have
\begin{eqnarray*}
\nb u_2^{(i)} &=& \sum_j \beta_{ij}\tensor u_1^{(j)} \\
\nb u_3 &=& \sum_i \gamma_i \tensor u_2^{(i)}
\end{eqnarray*}
for some $\beta_{ij}, \gamma_i \in \om$.
Differentiating the equation $\ag{u_1^{(i)},u_3} = 0$,
we get $\gamma_i = -d\log q_i$.

On the other hand,
since $\UU_2$ is a sub-crystal,
we can write
\[ \F u_2^{(i)} = p^2\left(a_iu_0 + \sum_j \left(b_{ij} u_1^{(j)}
			+ e_{ij} u_2^{(j)}\right)\right). \]
Here and in the rest of the proof,
we abbreviate $\F(\phi)\phi^*$ as $\F$.
Applying $\F$ to $\ag{u_1^{(i)},u_2^{(j)}} = \delta_{ij}$,
one gets $e_{ij} = \delta_{ij}$.

Suppose
\[ \F u_3 = p^3\left( cu_0 + \sum_i\left( f_iu_1^{(i)}
		+ e_iu_2^{(i)}\right) + fu_3\right). \]
Applying $\F$ to
$\ag{u_0,u_3} = 1, \ag{u_1^{(i)},u_3} = 0$ and $\ag{u_2^{(i)},u_3} = 0$,
one gets $f = 1, e_i = 1$ and $f_i = a_i$, respectively.

The relations \eqref{abck} between $a, b, c, \beta$
are derived from the flatness $\nb\F = \F\nb$.
Applying it to $u_3$ and $u_2^{(i)}$,
one obtains the first one and the last two equations, respectively.
\qed\medskip

\noin{\it Remark.}
For a discussion
about the possible relationship
between the constant term $c(0)$
of $c$ above
and a certain entry of the (complex) monodromy matrix
of the lifted algebraic family,
see \cite{Sh}
where the case of a family of quintic threefolds
has been computed explicitly.
The results obtained here do not provide any information of $c(0)$.

\subsection{Crystals and the mirror maps from degenerating families}

We now consider the situation
where the crystals are obtained from
certain families with a special kind of mild degenerations.
Given a flat projective scheme $X$ over $\Spec A$
with $X$ smooth over $W$,
we say that $X/A$ is a
{\em degenerating family of split type
(along the directions $t_i, 1\leq i\leq r$)}
if locally there exist elements
$\{x_{ij}\}_{1\leq i\leq r}^{1\leq j\leq n_i}$
in $\OO_X$ for some integers $n_i > 1$,
which, together with $p$,
form a part of a regular sequence,
such that $\OO_A \to \OO_X$ is given locally by
\[  t_i \mapsto \prod_{j=1}^{n_i} x_{ij} \quad(1\leq i\leq r). \]
Thus the fiber $Y$ over $t_1\cdots t_r = 0$
is a (reduced) normal crossing divisor of $X$,
and locally every irreducible component
of the fiber over any $t_i = 0, 1\leq i\leq m$,
is also of split type.

In this case, we equip $X$ and $A$
with the log structures
attached to the divisors $Y$ and $t_1\cdots t_r = 0$,
respectively.
Then $X/S$ is log-smooth;
the collection $\{x_{ij}^{-1}dx_{ij}\}$
form a part of a local basis
of the sheaf of log differential 1-forms on $X$
and we have
\begin{eqnarray*}
\om^1_S &\to& \om^1_X \\
\frac{dt_i}{t_i} &\mapsto& \sum_{j=1}^{n_i}\frac{dx_{ij}}{x_{ij}}
\quad(1\leq i\leq r).
\end{eqnarray*}

Let $X_0/S_0$ be the reduction mod $p$.
We have $H^\rho_{dR}(X/S) = H^\rho_{cris}(X_0/W)$,
where the cohomologies are understood to be the logarithmic ones.
Denote by $\om^i_{X/S}$
the sheaf on $X$ of relative log differential $i$-forms.
Now suppose, for all $i, j$ and $\rho$, that
$H^\rho_{dR}(X/S)$ and $H^j(X, \om^i_{X/S})$
are free $A$-modules
and the Hodge to de Rham spectral sequence
\[ E_1^{ij} = H^j(X, \om^i_{X/S}) \Longrightarrow H^{i+j}_{dR}(X/S) \]
degenerates at $E_1$ (see~\cite[Thm (4.12)(3)]{Kato_log}).
Then $H^\rho_{dR}(X/S)$
is a Hodge $F$-crystal of weight $\rho$ over $S_0$
provided that $p > \rho$
by the log analogue of a theorem of Mazur and Ogus
(\cite[Cor 8.3.3]{O} and cf.~\cite[Thm 8.26]{BO}).
\smallskip

Now suppose that
$K$ can be embedded into the field $\CC$ of complex numbers
and we fix one such embedding.
Assume that
\begin{itemize}
\item $r = m$,
\item the map \eqref{dag} is an isomorphism, and
\item the crystal $H_{dR}^{\rho}(X/S)$ is of Calabi-Yau type.
\end{itemize}
We define the {\em mirror maps}
$\tilde{q} = \{\tilde{q}_i\}_{i=1}^m$ as follows.

Denote by $\tilde{u}_0$
the image of $u_0 \in H_{dR}^{\rho}(X/S)$ in
\[ \HH_{\CC} = H_{dR}^\rho ((X/S)\tensor_W \CC)
	= H_{dR}^\rho (X/S)\tensor_W \CC. \]
Let $\tilde{u}_1^{(1)}, \cdots, \tilde{u}_1^{(m)}$
be elements in the Hodge filtration $\widetilde\Fil^1$ of $\HH_{\CC}$
satisfying
\begin{equation}\label{nb_u}
\nb(\tha_i)\tilde{u}_1^{(j)}\Big|_{t = 0}
	= \delta_{ij}\cdot \tilde{u}_0\Big|_{t = 0}
	\quad\text{for all $0\leq i,j\leq m$,}
\end{equation}
where $\delta =$ Kronecker's delta.
We define $\tilde{q}_i \in t_i(1+t_i\unit{\CC\dsb{t}})$
by requiring
\[ \nb \tilde{u}_1^{(i)} = d\log \tilde{q}_i \tensor \tilde{u}_0
	\quad (1\leq i\leq m). \]
We remark that
the elements $\tilde{q}_i$
can also be defined via the Picard-Fuchs differential equation
associated with $\HH_{\CC}$,
cf.~\cite[\S 1]{Y}.

\begin{Prop}\label{Prop:int_qc}
With notations and assumptions as above,
assume that the Hodge $F$-crystal $H_{dR}^{\rho}(X/S)$ is ordinary.
Via the inclusions $W \rc K \rc \CC$,
we have $\tilde{q}_i \in W\dsb{t}$ for all $i$.
\end{Prop}

\pf
We compare the $\tilde{q}_i$
with the $q_i \in W\dsb{t}$ defined in Thm \ref{Thm:weight=3}.
The condition \eqref{nb_u} forces that $\tilde{q}_i = a_iq_i$
for some constant $a_i \in \CC$.
Since $t_i^{-1}\tilde{q}_i(0) = 1$ by definition,
we have $a_i = (t_i^{-1}q_i)(0) \in \unit{W}$
and the assertion follows.
\qed

\subsection{Global case}
We now turn to the following global consideration.
Let $R$ be an open part
of the integral closure of a number field $R_\QQ$.
Suppose that we have a degenerating family $\pi: X \to \Spec R\dsb{t}$
of split type along $t_1,\cdots,t_m$
(defined in the similar way as in the beginning of \S \ref{Sect:App}$(b)$).
As before
we equip the family with the log structures
associated with the degenerate fiber of $X$
and the coordinate hyperplanes of the base $R\dsb{t}$,
respectively.
Let $\tilde{\om}$ denote
the sheaf on $X$ of the relative log differential 1-forms
of $X/R\dsb{t}$.
We assume that $H^j(X,\tilde{\om}^i)$ are free over $R\dsb{t}$
for all $i,j$, and
the relative Hodge to de Rham spectral sequence
degenerate at $E_1$.
These two requirements can be achieved by shrinking $R$.

Consider the relative de Rham cohomology
\[ \widetilde{\HH} = H^3_{dR}(X/R\dsb{t}) \]
of degree three
with the Hodge filtration $\widetilde{\Fil}^{\bullet}$
and the Gauss-Manin connection $\nb$.
Similar to \eqref{dag},
we have an $R\dsb{t}$-linear map
\begin{equation}
\gamma: \om_{R\dsb{t}}^{\vee} \to
	\Hom_{R\dsb{t}}(\widetilde\Fil^3, \widetilde\Fil^2/\widetilde\Fil^3),
\end{equation}
where $\om_{R\dsb{t}}$ denotes
the module of log differential 1-forms of $R\dsb{t}$.
We assume that
the map $\gamma$ is an isomorphism
and $\widetilde{\HH}$ is of Calabi-Yau type
(defined in a similar way as in the beginning of \S \ref{Sect:App}$(a)$).
Just as in the constructions of previous two subsections,
there exists a basis
$\{\tilde{u}_0, \tilde{u}_1^{(i)}, \tilde{u}_2^{(i)}, \tilde{u}_3\}_{i=1}^m$
of $\widetilde{\HH}$
and elements $\tilde{q} = \{\tilde{q}_i\}_{i=1}^m$
and $\{\tilde\kappa_{ijk}\}_{1\leq i,j,k\leq m}$ of $R\dsb{t}$
satisfying the following three conditions
(cf.~\cite[\S 14]{Del_inf} for the considerations
from the point of view of the Betti side):
\begin{enumerate}
\item
$\tilde{u}_0$ generates the unique rank one $R$-module
of horizontal sections in $\widetilde{\HH}$.
\item
We have $\tilde{q}_i \in t_i(1+t_i\unit{R\dsb{t}})$
and $\nb\tilde{u}_1^{(i)} = d\log \tilde{q}_i \tensor \tilde{u}_0$
for $1\leq i\leq m$
as in the previous subsection.
\item
Together with $\{\tilde{u}_0,\tilde{u}_1^{(i)}\}_{i=1}^m$,
the reminder $\{\tilde{u}_2^{(i)},\tilde{u}_3\}_{i=1}^m$
form a symplectic basis of $\widetilde\HH$
similar to the description in Thm \ref{Thm:weight=3}.
Define similarly
\begin{equation*}
\tilde{\kappa}_{ijl} \in R\dsb{t} = R\dsb{\tilde{q}}
\end{equation*}
corresponding to the elements $\kappa_{ijl}$
given in \eqref{eqn:def:kappa}.
\end{enumerate}

For a prime $\wp$ of $R$,
denote by $R_{\wp}$ its completion at $\wp$
and by $X_{\wp}/R_{\wp}\dsb{t}$ the base change of $X/R\dsb{t}$.
Consider the crystal
\[ \HH_\wp = H^3_{dR}(X_\wp/R_\wp\dsb{t})
	= \widetilde\HH \tensor_{R} R_\wp. \]
Let $e: R_\wp\dsb{t} \to R_\wp$ be the augmentation.
Denote by $\bar{u}_0 \in e^*\HH_\wp$
the induced element of $\tilde{u}_0$
through the base changes.
As an application of previous results,
we have the following.

\begin{Thm}\label{Thm:kappa=kappa'}
With notations and assumptions as above,
suppose that for each prime $\wp$ of $R$,
the crystal $\HH_\wp$ is ordinary
and $\bar{u}_0$ and $(t_i^{-1}q_i)(0), 1\leq i\leq m$,
are fixed by the absolute Frobenius,
where $q_i$ are defined in Theorem \ref{Thm:weight=3}.
Consider the expansion
\[ \tilde{\kappa}_{ijl} = \eps_{ijl}(0) + \sum_{|n| > 0} \eps_{ijl}(n)
	\frac{n_i n_j n_l \tilde{q}^n}{1-\tilde{q}^n}
		\quad(1\leq i,j,l\leq m, \eps_{ijl}(n) \in R_\QQ). \]
Then $\eps_{ijl}(0) \in R \dc \QQ$
and $\tr (\eps_{ijl}(n)) \in \tr R$
for all $n\neq 0$,
where $\tr : R_\QQ \to \QQ$ is the trace map.\footnote{
The $\eps_{ijl}(n)$ are the {\it instanton numbers}.}
\end{Thm}

\pf
Indeed, on $\HH_\wp$,
the two bases
$\{\tilde{u}_0, \tilde{u}_1^{(i)}, \tilde{u}_2^{(i)}, \tilde{u}_3\}$,
induced from that of $\widetilde{\HH}$
by the natural base change,
and
$\{u_0, u_1^{(i)}, u_2^{(i)},u_3\}$,
constructed in Thm \ref{Thm:weight=3},
only differ by a multiplicative $p$-adic integer.
More precisely,
let $W$ be the ring of Witt vectors
of an algebraic closure of $R/\wp$.
Then there exists a non-zero $\alpha \in W$
such that $u_0 = \alpha\tilde{u}_0$.
Since $\F\bar{u}_0 = \bar{u}_0$,
we obtain $\alpha \in \ZZ_p$.
Consequently
we have $u_1^{(i)} = \alpha\tilde{u}_1^{(i)}$
and $u_2^{(i)} = \alpha^{-1}\tilde{u}_2^{(i)}$ for all $i$.
Thus $\kappa_{ijl} = \alpha^2\tilde{\kappa}_{ijl}$.
On the other hand,
we have $R_\wp\dsb{\tilde{q}} = R_\wp\dsb{q}$
and $\tilde{q_i}\frac{\partial}{\partial\tilde{q}_i} = q_i\frac{\partial}{\partial q_i}$
by the proof of Prop \ref{Prop:int_qc}.
The assertion then follows from
Thm \ref{Thm:weight=3} above
and Lemma \ref{Lemma:KSV} below.
(Note that if the Frobenius sends $q_i$ to $q_i^p$,
then it sends $\tilde{q}_i$ to $\tilde{q}_i^p$
by our assumption on the leading terms of $q_i$.)
\qed

\begin{Lemma}[{\cite[Lemma 1, 2]{KSV}}]\label{Lemma:KSV}
Let $V$ be the ring of Witt vectors
of a finite field of characteristic $p$
and $\tr: \QQ\tensor_{\ZZ}V \to \QQ_p$ the trace map.
Fix $1 \leq i_j \leq m$ for $j=1,\cdots,r$.
Let $\kappa \in V\dsb{q}$ be written formally as
\[ \kappa = \eps_0 + \sum_{|n| > 0} \eps_n
	\frac{n_{i_1}n_{i_2}\cdots n_{i_r} q^n}{1-q^n}
	\quad (n = (n_l)_{l=1}^m \in \NN^m, \eps_n \in \QQ\tensor_{\ZZ}V). \]
Let $\phi$ be the lifting of the absolute Frobenius on $V\dsb{q}$
given by $\phi(q) = q^p$
and $\pt_i = q_i\frac{\pt}{\pt q_i}$.
Suppose
\begin{equation}\label{kappa-f}
\phi(\kappa) - \kappa = \pt_{i_1}\pt_{i_2}\cdots\pt_{i_r} f
	\quad\text{for some $f \in V\dsb{q}$}.
\end{equation}
Then
\[ \tr \eps_n \in \ZZ_p
	\quad\text{for all $n \neq 0$}. \]
\end{Lemma}

\pf
Let $s$ be the rank of $V$ over $\ZZ_p$.
The assertion in the case $(r,s) = (3,1)$
is proved by using the M\"obius inversion formula
in \cite[\S 3]{KSV} (see especially Equation (30) therein).
The general case is proved similarly.
We remark that in case $s = 1$,
the condition \eqref{kappa-f} is also necessary.
\qed\medskip

\noin{\em Remarks.}
(i) One know that
if the irreducible components of the central fiber of the family
$X_\wp/R_\wp\dsb{t}$
are ordinary,
the crystals $H^\rho_{dR}(X_\wp/R_\wp\dsb{t})$ are ordinary
(\cite[3.23]{Mo}).
With notations as above,
it is easily seen that
the cohomology group $\widetilde{\HH}$ is of Hodge-Tate type
in the sense of \cite[\S 6]{Del_inf}
if $\HH_\wp$ is ordinary for some $\wp$ by Prop \ref{Prop:UH}.

(ii) The conditions on $\bar{u}_0$ and $(t_i^{-1}q_i)(0)$
in Thm \ref{Thm:kappa=kappa'}
can be checked by looking at the $F$-crystal $e^*\HH_\wp$.
On one hand,
$\F \bar{u}_0$ is a constant multiple of $\bar{u}_0$;
on the other hand,
the leading terms of $q_i$ are normalized by \eqref{club}
(see also \eqref{bigstar}).
However, the log structure on $e^*X$ depends on
its embedding to $X$.

\appendix
\section{Appendix: The group structure}

In the classical case,
the deformation functor to local $W$-algebras
of an ordinary abelian variety or K3 surface $X$
over the field $k$
has a natural formal group structure,
and this leads to the notion of the canonical lifting of $X$
(see \cite{Del_CC}).
In this appendix,
we provide an analogous group structure
on the deformation of an ordinary crystal of weight one
under the appearance of a logarithmic structure.
The explanation here is not satisfactory yet
since the group functor is not (pro-)representable in this case.
\smallskip

Fix $r \in \NN$.
In this section,
let $\{e_i\}_{i=1}^r$ be the standard basis of $\NN^r$.

Let $\widetilde{\Art}$
be the category of
Artinian local $\ZZ_p$-algebras $(R,\mathfrak{m})$
with fixed inclusions $R/\mathfrak{m} \hookrightarrow k$,
equipped with log structures $\MM_R$
and structure morphims sitting in the sequence
\begin{equation}\label{A:str-morph}
(\ZZ_p,\unit{\ZZ_p}) \to (R,\MM_R) \to \logpt{k}{r}
\end{equation}
of log rings.
The morphisms of $\widetilde{\Art}$
are the morphisms of log rings
compatible with the structure morphisms \eqref{A:str-morph}.

For an Artinian local $W$-algebra $(R,\mathfrak{m})$
with a fixed $R/\mathfrak{m} \to k$,
denote by $\logpt{R}{r}$ the ring $R$ equipped
with the log structure
\begin{eqnarray}\label{A:dir-sum}
\nnunit{r}{R} &\to& R \\
\nonumber (n,x) &\mapsto& 0^n\cdot x.
\end{eqnarray}
and with the structure morphism
$({\rm id,mod}\, \mathfrak{m}): \nnunit{r}{R} \to \nnunit{r}{k}$.
(We also extend this symbol to
pro-Artinian local $W$-algebras.)
Let $\Art$ be the subcategory of $\widetilde{\Art}$
consisting of such Artinian $\logpt{R}{r}$'s
whose morphisms consist of those
which preserve the direct sums \eqref{A:dir-sum} of the log structures.
Notice that $\Art$ does not have an initial object.
However each $\logpt{R}{r}$
admits a morphism $\logpt{(\ZZ_p)}{r} \to \logpt{R}{r}$
whose log part $\nnunit{r}{\ZZ_p} \to \nnunit{r}{R}$
is the direct sum of the identity and the structure morphism.
We denote this special object (with the arrows)
by $\mathring{\ZZ}_p$.

\subsection{The group functor $\SS$}

We define functors $\GG, \TT$ and $\SS$ on $\Art$,
which are pro-representible in $\widetilde{\Art}$.
The functor $\SS$ will be an extension of $\TT$ by $\GG$
and will provide the group structure
attached to a certain ordinary crystal
on the log scheme $S=(\Spec A, \LL)$
discussed in the following subsection.\smallskip

Let $\GG = \Spf \ZZ_p\dsb{z}, z = \{z_i\}_{i=1}^r$,
which is equipped with a binary operator $*:\GG\times\GG \to \GG$;
here
the log structure on $\GG$ is given by
\begin{eqnarray*}
\nnunit{r}{\ZZ_p\dsb{z}} &\to& \ZZ_p\dsb{z} \\
	(n,f) &\mapsto& z^nf
\end{eqnarray*}
with the structure morphism
$({\rm id,mod}\, (p,z)): \nnunit{r}{\ZZ_p\dsb{z}} \to \nnunit{r}{k}$,
and $*$ is induced
by the diagonal map on the level of rings and monoids.\footnote{
Notice however that $*$ is not representable.}
The $*$ is associative and commutative.
Moreover $\GG$ is a group functor.
Indeed
let $\mathring{e}$ be the unique element in
$\GG(\mathring{\ZZ}_p)$.
The pair $(*,\mathring{e})$ then turns $\GG$
into a monoid with unit.
Explicitly,
for $(\logpt{R}{r},\mathfrak{m})$ an object of $\Art$,
we have
\begin{equation}\label{eqn:S(R)}
\GG(\logpt{R}{r},\mathfrak{m}) = \left\{
	\begin{array}{rcl|c}
	\nnunit{r}{\ZZ_p\dsb{z}} &\to& \nnunit{r}{R} & \alpha_i \equiv 1 \\
	e_i &\mapsto& (e_i,\alpha_i) & \mod{\mathfrak{m}}
	\end{array} \right\};
\end{equation}
here the condition $\alpha \equiv 1 \mod{\mathfrak{m}}$
follows from the commutativity of the diagram
\[ \xymatrix{
(\ZZ_p\dsb{z}, \nnunit{r}{\ZZ_p\dsb{z}}) \ar[r]\ar[dr] & \logpt{R}{r} \ar[d] \\
	& \logpt{k}{r}. }\]
We identify $\GG(\logpt{R}{r})$ with the set $(1+\mathfrak{m})^r$
via the $r$-tuples $(\alpha_i)$ appeared in \eqref{eqn:S(R)}.
Under this identification,
$*$ corresponds the standard multiplication on $(1+\mathfrak{m})^r$.
We remark that
the map $\Spf \logpt{(\ZZ_p)}{r} \to \GG$,
whose log part reads
$({\rm id,mod}\, (z)): \nnunit{r}{\ZZ_p\dsb{z}} \to \nnunit{r}{\ZZ_p}$,
induces an isomorphism as functors on $\Art$.
\smallskip

Let $\TT = \Spf \ZZ_p\dsb{x -1}$, $x = \{x_i\}_{i=r+1}^m$,
be the formal torus over $\Spf \ZZ_p$ of dimension $(m-r)$,
with characters generated freely by $x$.
We regard $\TT$ as a group functor on $\Art$
by sending $\logpt{R}{r}$ to $\TT(R)$.
The assignment
\begin{eqnarray}\label{eqn:changebasis-T}
{\rm GL}_{m-r}(\ZZ_p) &\to& \{\text{basis of the characters of $\TT$}\}
	\nonumber \\
(\gamma_{ij})_{r< i,j\leq m} &\mapsto&
	\left\{\prod_{j=r+1}^m x_j^{\gamma_{ij}}\right\}_{i=r+1}^m
\end{eqnarray}
establishes a bijection.\smallskip

Let $B$ be a $\ZZ_p\dsb{x -1}$-algebra
where 
\begin{equation}\label{eqn:S=GT}
B \cong \ZZ_p\dsb{x -1} \widehat{\tensor}_{\ZZ_p} \ZZ_p\dsb{z}
	\quad\text{non-canonically}
\end{equation}
and is equipped with the log structure
attached to the divisor $z=0$.
Let $\SS = \Spf B$
regarded as a group functor
formally sitting in the short exact sequence
\[\xymatrix{
1 \ar[r] & \GG \ar[r] & \SS \ar[r] & \TT \ar[r] & 1. }\]
Now fix a collection $\tilde{z} = \{\tilde{z}_i \in B\}_{i=1}^r$ of liftings of $z_i$.
By an inspection of the group $\SS(\logpt{(\ZZ_p)}{r})$,
we see that
different choices of the splitting \eqref{eqn:S=GT} of $B$
are in one-to-one correspondence
with the changes of the collection $\tilde{z}$ in the way
\begin{equation}\label{eqn:changebasis-G}
\tilde{z}_i \mapsto x^{\gamma_i} \tilde{z}_i,
	\quad\text{for some $\gamma_i \in (\ZZ_p)^{m-r}, 1\leq i\leq r$}.
\end{equation}

\subsection{The structure of $S$ attached to a crystal}

We will consider the following situation
of an ordinary crystal with
an explicitly described connection and Frobenius structure.
In the following,
we change the notations of the variables
used in the main content slightly.

Let $\HH$ be an ordinary Hodge $F$-crystal of weight one over $S_0$
with an $A$-basis $\{a_i,b_j\}$,
where $(1,1)\leq(i,j)\leq(g,h)$ with $gh = m$, satisfying
\begin{eqnarray}\label{eqn:App-H-Fil}
\Fil^1 &=& \text{the $A$-module generated by $\{b_j\}_{j=1}^h$} \\
\nb a_i = 0 &;& \nb b_j = \sum_{i=1}^g d\log q_{ji} \tensor a_i
	\quad\text{for some $q_{ji} \in A$ with $d\log q_{ji} \in \om$.}
\end{eqnarray}
Suppose that there exists a subset $E$ with $r$ elements
such that
\begin{eqnarray}
A = W\dsb{q_{ji} - \eps_{ji}} &;& \LL = \bigcup_{n \in \NN^E} q^n \unit{A}
\end{eqnarray}
where $q = (q_{ji})$ and
\[ \eps_{ji} = \left\{\begin{array}{cl}
	0 & \text{if $(i,j) \in E$} \\
	1 & \text{otherwise}. \end{array}\right. \]
Denote by $\phi$ the lifting to $S$ of the absolute Frobenius
given by $\phi(q) = q^p$.
Further suppose that the Frobenius structure of $\HH$ is given by
\begin{eqnarray}\label{eqn:App-H-F}
\F(\phi)\phi^*a = a &;& \F(\phi)\phi^*b_j = pb_j.
\end{eqnarray}
Now fix a labeling on the elements of $E$ and its complement
by $\{1,\cdots,r\}$ and $\{r+1,\cdots,m\}$, respectively.

\begin{Prop}[{Cf.~\cite[Th 2.1.14]{Del_CC}}]
Attached to the ordinary crystal $\HH$ of weight one as above,
one has a canonical isomorphism of log schemes
\begin{equation}\label{eqn:S-SS}
S \cong \SS,
\end{equation}
where $\SS$ is defined in \eqref{eqn:S=GT}.
\end{Prop}

\pf
With notations as above,
the desired isomorphism $\SS \to S$ is given by
\begin{equation*}
q_l \mapsto \left\{\begin{array}{cl}
	z_l & \text{if $l\leq r$} \\
	x_i & \text{if $l>r$}. \end{array}\right.
\end{equation*}
Notice that different choices of the possible $\{a,b\}$ and $q$'s for $\HH$
satisfying the conditions \eqref{eqn:App-H-Fil} - \eqref{eqn:App-H-F}
are exactly in one-to-one correspondence
with the combinations of the transformations
\eqref{eqn:changebasis-T} and \eqref{eqn:changebasis-G} for $\SS$.
Thus the isomorphism \eqref{eqn:S-SS}
is indeed canonical
(cf.~\cite[Th 2.1.14]{Del_CC}).
\qed

\subsection{The periods of the $\logpt{W}{r}$-points of $S$}

With notations and assumptions as in the previous subsection,
let $\overline{\HH}_0 =$ the $F$-crystal over $\logpt{k}{r}$
obtained by restriction of $\HH$ to the origin
(and forgetting the Hodge filtration).
Then $\overline{\HH}_0$ has the basis
$\{\bar{a},\bar{b}\} :=$ the restriction of $\{a,b\}$ to 0,
with
\begin{eqnarray*}
\F \bar{a}_i = \bar{a}_i &;& \F \bar{b}_j = p\bar{b}_j
	\quad((1,1)\leq(i,j)\leq(g,h)) \\
\N_{ls} \bar{a}_i = 0 &;& \N_{ls} \bar{b}_j = \delta_{lj} \bar{a}_s
	\quad((s,l) \in E).
\end{eqnarray*}
Let $e_\phi$ be the \Teich lifting attached to $\phi$
and $\Fil_\phi^\bullet$ the Hodge filtration
of the ordinary crystals $e_\phi^*\HH$ over $\logpt{k}{r}$.
The $\bar{b}$ forms a basis of $\Fil_\phi^1$ over $W$.

Let $\chi \in S(\logpt{W}{r})$
be a $\logpt{W}{r}$-point of $S$
and $\Fil_\chi^\bullet$ the Hodge filtration of $\chi^*\HH$.
Explicitly with the notations in the proof of Lemma \ref{Lemma:Teich-Base}, the element
\[ \chi = (\chi^{\flat},\chi^{\sharp}): (A,\LL) \to \logpt{W}{r} \]
is determined by specifying $\alpha_{ji} \in pW$ and  $\beta_{ji} \in 1+pW$
such that
\begin{eqnarray*}
\chi^\flat(q_{ji}-1) &=& \alpha_{ji} \quad \text{for $(i,j) \not\in E$} \\
\chi^\sharp(e_{ji},1) &=& (e_{ji},\beta_{ji})
	\quad \text{for $(i,j) \in E$}.
\end{eqnarray*}
In order to put all the $q$'s into the same framework,
we write
\[ \beta_{ji} := \chi^\flat(q_{ji}) \in 1+pW
	\quad\text{for $(i,j) \not\in E$}. \]
Thus we have an isomorphism
\begin{eqnarray*}
S(\logpt{W}{r}) &\cong& (1+pW)^m \\
\chi &\mapsto& (\beta_{ji}).
\end{eqnarray*}
Under this,
the group structure, via the identification \eqref{eqn:S-SS}, of the left side
coincides with the component-wise multiplication on the right side.

\begin{Prop}\label{Prop:period}
With notations as above
and under the canonical identification
\[ \chi^*\HH \cong e_\phi^*\HH = \overline{\HH}_0, \]
given in \eqref{eqn:e^*H},
the $W$-module $\Fil_\chi^1$
(regarded as a submodule of $\overline{\HH}_0$)
is generated by the elements
\[ \bar{b}_j + \sum_{i=1}^g \varpi_{ji}\bar{a}_i
	\quad (1\leq j\leq h) \]
where $\varpi_{ji} = \log \beta_{ji} \in pW$.
\end{Prop}

\pf
The $W$-submodule $\Fil_\chi^1$ of $\chi^*\HH$ is generated by
$\chi^*b$.
Inductively on $n$, we obtain
\begin{eqnarray*}
\prod_{s=0}^{n-1}\left(\nb\left(q_{ji}\frac{d}{dq_{ji}}\right)-s\right)b_l
	&=& q_{ji}^n\nb\left(\frac{d}{dq_{ji}}\right)^nb_l \\
	&=& (-1)^{n-1}(n-1)!\delta_{jl}a_i
	\quad (n > 0, 1\leq l\leq h).
\end{eqnarray*}
Thus by the explicit formula \eqref{eqn:e^*H} for the identification,
the assertion follows.
\qed\medskip

We call $\varpi_{ji}$ the {\it periods}
associated with $\chi$.
Let $\chi' \in S(\logpt{W}{r})$ be another morphism
corresponding to $(\beta_{ji}')$ with periods $\varpi_{ji}'$
as above.
Then clearly the sum $\chi+\chi' \in S(\logpt{W}{r})$,
via the group structure \eqref{eqn:S-SS},
has associated periods $\varpi_{ji}+\varpi_{ji}'$.
One checks readily that
the neutral element $\chi_0$ is characterized
by saying that the induced Hodge
and the slope filtrations on $\chi_0^*\HH$ coincide.
This is the log analogue of \cite[2.1.3, 2.1.8]{Del_CC}.\medskip

\noin{\it Remark.}
The category $\Art$ considered in this section
does not contain the interesting log ring
$W' := (W$ with the log structure attached to the divisor $\Spec k$).
It is tempting to find a better explanation of this exclusion
or find a way to incorporate the direction ``$p$"
into the deformation.
Notice however that
the Frobenius on $\logpt{k}{1}$
does not admit a lifting to $W'$.

\end{document}